\documentclass[14pt,a4paper,openright,reqno]{extarticle}
\usepackage[english]{babel}
\usepackage{newlfont}
\usepackage{amsbsy}
\usepackage[dvips]{graphicx}
\usepackage{color}
\usepackage{bbm}
\usepackage{graphicx, subfigure}
\usepackage{color}
\usepackage[center]{caption2}
\usepackage{amssymb}
\usepackage{amsmath}

\newcommand{\rk}{{\operatorname{rk}}}

\newcommand{\mi}{{\operatorname{mid }}}
\newcommand{\rad}{{\operatorname{rad }}}

\usepackage{latexsym}
\usepackage{amsthm}
\usepackage{amsthm}

\usepackage[dvips]{geometry}
\theoremstyle{plain}
\newtheorem{thm}{Theorem}
\newtheorem{cor}[thm]{Corollary}

\newtheorem{lem}[thm]{Lemma}

\newtheorem{nota}[thm]{Notation}

\theoremstyle{remark}
\newtheorem{rem}[thm]{Remark}

\newcommand{\R}{\mathbb{R}}
\newcommand{\Q}{\mathbb{Q}}

\newcommand{\N}{\mathbb{N}}

\newcommand{\I}{\mathbb{I}}
\newcommand{\al}{\pmb{\alpha}}

\usepackage{mathtools}
\def\alphaltiset#1#2{\ensuremath{\left(\kern-.2em\left(\genfrac{}{}{0pt}{}{#1}{#2}\right)\kern-.2em\right)}}
\usepackage[center]{caption2}

\begin{document}

\title{Interval matrices: realization of ranks by rational  matrices }
\author{Elena Rubei}
\date{}
\maketitle

{\footnotesize\em Dipartimento di Matematica e Informatica ``U. Dini'', 
viale Morgagni 67/A,
50134  Firenze, Italia }

{\footnotesize\em
E-mail address: elena.rubei@unifi.it}

\def\thefootnote{}
\footnotetext{ \hspace*{-0.36cm}


{\bf \copyright } This manuscript version is made available under the CC-BY-NC-ND 4.0 license https://creativecommons.org/licenses/by-nc-nd/4.0/}

\begin{abstract} 
Let $\al $ be  a $p \times q$ interval matrix with $p \geq q$ 
 and   with  the endpoints of all its entries  in $\Q$.  We prove that, if  $\al $ contains a rank-$r$ real matrix 
 with $r \in \{2, q-2,q-1,q\}$, then it contains   a rank-$r$ rational matrix.
\end{abstract}

\section{Introduction}

Let $p , q \in \N \setminus \{0\}$;
a    $p \times q$  interval matrix $\al$ is a  $p \times q$  matrix 
whose entries are intervals in $\R$; 
we usually denote the entry $i,j$, $\al_{i,j}$, by $[\underline{\alpha}_{i,j}, \overline{\alpha}_{i,j}]$ with $\underline{\alpha}_{i,j} \leq 
\overline{\alpha}_{i,j}$ and we point out that we denote every interval matrix in bold.
A  $p \times q$  matrix $A$  with entries in $\R$ is  said contained in  a $p \times q$ interval matrix $  \al $ if $a_{i,j} \in \al_{i,j} $ for any $i,j$. 
There is a wide literature about interval matrices 
and the rank of the matrices they contain.
In this paper we consider the following problem:
let $\al $ be  an interval matrix  whose entries have rational endpoints; 
for which $r$ can we deduce that, if $\al$ contains a rank-$r$ real matrix, then  $\al$ contains a rank-$r$ rational matrix?
 
Before sketching our results, we illustrate shortly some of the literature 
on interval matrices, and the rank of the contained matrices,  on  partial matrices and on the  matrices with a given sign pattern; these last two research fields are connected with the theory of interval matrices. 
 
  Two
of the most famous theorems on interval matrices are Rohn's theorems
on full-rank interval matrices.
We say that a  $p \times q $ interval matrix $\al$ has full rank if and only if  all the matrices contained in $\al$ have rank equal to $\min\{p,q\}$. 
For any  $p \times q $ interval matrix $\al = ( [\underline{\alpha}_{i,j}, \overline{\alpha}_{i,j}])_{i,j}$ with $\underline{\alpha}_{i,j} \leq \overline{\alpha}_{i,j}$,
let $\mi(\al)$, $\rad(\al)$ and $|\al|$ be respectively the midpoint, the radius and the modulus of $\al$, that is 
 the  $p \times q$   matrices such that $$ \mi (\al)_{i,j}= \frac{\underline{\alpha}_{i,j}+ \overline{\alpha}_{i,j}}{2}, \hspace*{1cm} \rad(\al)_{i,j}= \frac{\overline{\alpha}_{i,j}- \underline{\alpha}_{i,j}}{2}, $$ 
$$|\al|_{i,j} = \max\{|\underline{\alpha}_{i,j}|,| \overline{\alpha}_{i,j}|\} $$ for any $i,j$. 
The following theorems
characterize respectively full-rank {\em square} interval matrices and 
full-rank $p \times q$ interval matrices, 
 see  \cite{Rohn},
\cite{Rohn3}, \cite{Rohn4}, \cite{Shary}; 
see  \cite{Rohn} and  \cite{Rohn2} for other characterizations. 

\begin{thm} {\bf (Rohn)} \label{Rohn1}
Let   $\al =( [\underline{\alpha}_{i,j}, \overline{\alpha}_{i,j}])_{i,j}$ be  a  $p \times p$  interval matrix, where
 $\underline{\alpha}_{i,j} \leq \overline{\alpha}_{i,j}$  for any $i,j$.
Let $Y_p=\{-1,1\}^p$ and, for any $x
 \in Y_p$, denote by $T_x$ the diagonal matrix whose diagonal is $x$.
Then  $\al$ is a full-rank interval matrix 
if and only if,  for each
$x,y \in Y_p$,  $$\det\Big(\mi(\al)\Big) \,\det\Big(\mi(\al) - T_x \, \rad(\al) \, T_y\Big)>0. $$
\end{thm}

\begin{thm}  {\bf (Rohn)} 
A $p \times q$ interval matrix $\al$ with $p \geq q$ has full rank if and only if the system
of inequalities $$\hspace*{2cm} |\mi(\al) \, x| \leq \rad(\al) \,
|x|, \hspace*{1.5cm} x \in \R^q$$ has only the trivial solution $x=0$. 
\end{thm}

  A research area which can be connected with the theory of interval matrices 
  is the one of the {\em partial matrices}: let $K$ be a field; a partial matrix over $K$ is a 
matrix where only some of the entries are given and they are elements of 
$K$;
a completion of a partial  matrix is a specification of the unspecified entries. 
In  \cite{CJRW}, Cohen,  Johnson,  Rodman and Woerdeman determined 
the maximal rank of the completions of a partial
matrix in terms of the ranks and the sizes of its maximal  
specified submatrices; see also \cite{CD}
for the proof.
 The problem of a theoretical characterization of   the minimal rank of   the completions of a partial matrix seems more difficult  and it has been solved only in some particular cases. 
 We quote also the papers \cite{Peeters} and
 \cite{Shi2} about the NP-hardness of the problem and the paper \cite{HHW} for rank-$1$ completions.
 
  
  In \cite{Ru2} we generalized Theorem \ref{Rohn1} to  matrices whose entries are closed connected nonempty subsets of $\R$, i.e. the so-called matrices in Kahan arithmetic.
  
In \cite{Ru1}  we 
 determined the maximum rank of the matrices contained in a given interval matrix
and we gave  a theoretical
      characterization
 of interval matrices containing at least a  matrix of rank $1$. 
In the previous paper \cite{G-S}, the authors studied the complexity of an algorithm to decide if an interval matrix contains a rank-one matrix and proved that the problem is NP-complete.

Finally we quote another research area which can be related 
to partial matrices, to interval matrices and, more generally,
to general interval matrices: the one of the {\em matrices with a given sign pattern}; let $Q$ be a $p \times q$ matrix with entries in $\{+,-,0\}$; we say that $A \in M(p \times q, \R)$ has sign pattern $Q$ if, for any $i,j$, we have that  $a_{i,j} $ is 
positive (respectively negative, zero) if and only if $Q_{i,j}$ is 
$+$ (respectively $-,0$). Obviously the set of the matrices 
with a given sign pattern can be thought as a matrix whose entries are in  $\{ (0, +\infty), (-\infty,0) , [0]\}$.
There are several papers 
studying the minimal and maximal rank of the matrices with a given sign pattern, see for instance \cite{A1}, \cite{A2}, \cite{A3},  \cite{Shi}. In particular, in \cite{A1} and \cite{A2}
the authors proved that the minimum rank of the real matrices with a given sign pattern is realizable by a rational matrix
in case this minimum is at most $2$ or at least $\min\{p,q\}-2$.


Obviously the three theories we have quoted, that is 
the theory of interval matrices, the theory of partial matrices,  and the theory of matrices with a given sign pattern can be seen as parts of the same theory: the one of {\em subset matrices}, i.e. matrices whose entries are subsets of a given field; we denote also subset matrices in bold.

In \cite{Ru3} we proved the following theorems:

\begin{thm} \label{Q-R1}
Let  $ p \geq q$ and let   
$\al =( [\underline{\alpha}_{i,j}, \overline{\alpha}_{i,j}])_{i,j}$ be a $p \times q$ interval matrix with 
$\underline{\alpha}_{i,j} \leq \overline{\alpha}_{i,j}$ and 
 $\underline{\alpha}_{i,j}, \overline{\alpha}_{i,j} \in \Q$ for any $i,j$. 
If there exists $ A \in \al$ with $\rk(A) <q$, then 
there exists $ B \in \al \cap M(p \times q, \Q)$ with $\rk(B) <q$.
\end{thm}

\begin{thm} \label{Q-R2}
Let  $ p \geq q$ and let   
$\al =( [\underline{\alpha}_{i,j}, \overline{\alpha}_{i,j}])_{i,j}$ be a $p \times q$ interval matrix with 
$\underline{\alpha}_{i,j} \leq  \overline{\alpha}_{i,j} $
and $\underline{\alpha}_{i,j}, \overline{\alpha}_{i,j} \in \Q$ for any $i,j$. 
If there exists $ A \in \al$ with $\rk(A) =1$, then 
there exists $ B \in \al \cap M(p \times q, \Q)$ with $\rk(B) =1$.
\end{thm}

 Moreover,  in \cite{Ru3} we observed (see in Remark 13 there) that from the papers \cite{Ber}, \cite{Shi} and \cite{Ko} we can deduce that
    it is not true that, for any $r$,
 if a $p \times q$ interval matrix   with  the endpoints of all its entries  in $\Q$  contains a rank-$r$ real matrix, then it contains   a rank-$r$ rational matrix. In particular this is not true for $r=3, \min\{p,q\}-3$.
 
 In this paper we prove that,
 if a $p \times q$ interval matrix with $p \geq q$ 
 and   with  the endpoints of all its entries  in $\Q$  contains a rank-$r$ real matrix, then it contains   a rank-$r$ rational matrix 
 for $r =2,q-2,  q-1,q$, see Theorem \ref{Q-R2}, Theorem 
 \ref{Q-Rq-1} and Remark \ref{Q-Rq}.
 Summarizing we get the following result; observe that the behaviour of interval matrices is similar to the one of the matrices with a given sign pattern
 showed in \cite{A1} and \cite{A2},
  even if, to prove it, we have to use a  technique which is 
  different from the one in \cite{A1} and \cite{A2}.

\begin{thm} \label{summarize}
Let  $ p \geq q$ and let   
$\al =( [\underline{\alpha}_{i,j}, \overline{\alpha}_{i,j}])_{i,j}$ be a $p \times q$ interval matrix with 
$\underline{\alpha}_{i,j} \leq  \overline{\alpha}_{i,j} $
and $\underline{\alpha}_{i,j}, \overline{\alpha}_{i,j} \in \Q$ for any $i,j$. 
If $r \in \{0,1,2,q-2,q-1,q\}$ and  there exists $ A \in \al$ with $\rk(A) =r$, then 
there exists $ B \in \al \cap M(p \times q, \Q)$ with $\rk(B) =r$.
\end{thm}
 
 
 
\section{Notation and first remarks}

$\bullet$  Let $\R_{>0}$ be the set $\{x \in \R | \; x >0\}$ and
let $\R_{\geq 0}$ be the set $\{x \in \R | \; x  \geq 0\}$; we define analogously $\R_{<0}$ and $\R_{ \leq 0}$. We denote by $\I$ the set $\R -\Q$. 

$\bullet $ Throughout the paper  let $p , q \in \N \setminus \{0\}$.

$\bullet $  For any set $X$, let $|X|$ be the cardinality of $X$.

$\bullet$ For any field $K$, let $M(p \times q, K) $ denote the set of the $p \times q $  matrices with entries in $K$. For any $A \in M(p \times q, K)$, let $\rk(A) $ be the rank of $A$,
 let $A^{(j)}$ be the $j$-th column of $A$ and, more generally, let $A^{(j_1,\ldots,j_r)}_{(i_1, \ldots, i_s)}$ be the submatrix of $A$ given by the columns $j_1,\ldots, j_r$ and the rows $i_1, \ldots, i_s$ of $A$ with the orders,  respectively, $j_1,\ldots, j_r$ and $i_1, \ldots, i_s$.

$\bullet $
For any vector space $V$ over a field $K$  and any $v_1,\dots, v_k \in V$, let $\langle v_1,\dots, v_k\rangle $ denote the 
span of $v_1,\dots,v_k$.

$\bullet $
 Let   $\al $ be a $p \times q $  subset matrix
 over a field $K$.  
 Given a matrix $A \in M(p \times q, K) $, we say that  $ A \in \al $ if and only if  $a_{i,j} \in \al_{i,j} $ for any $i,j$.
  
We say that an  entry of $\al$ is  {\bf degenerate} if its cardinality  is $1$.

$\bullet$ 
 Let   $\al $ and $\al'$ be two $p \times q $  interval matrices.  We say that $$ \al' \subset \al$$ if 
 $ \al'_{i,j} \subset \al_{i,j}$ for every $i,j$.

We defer to some classical books on interval analysis, such as \cite{Moore}, \cite{Neu} and \cite{Moore2} for the definition of sum and multiplication of two intervals. In particular, 
for any interval $\alpha$ in $\R$ and any interval $\beta$ either in 
$ \R_{>0}$ or in $\R_{<0}$,
we define $\frac{\alpha}{\beta}$ to be the set $\left\{\frac{a}{b} | \; 
a \in \alpha, \; b \in \beta \right\}$.

\section{Rational realization of the rank $2$}

\begin{lem}
Let $K$ be a field and let $k,n \in \N\setminus \{0\}$. Let $A \in M(k \times n , K)$ with $ n > k$. If $A^{(1,\ldots, k)}$ is invertible, 
then a basis of
the kernel of $A$ is given by the following vectors
in $K^n$ for $j=k+1, \ldots, n$:
$$v_j:= \begin{pmatrix} 
\det (A^{(2, \ldots , k, j)}) \\
- \det ( A^{(1,3, \ldots , k, j)}) \\
\vdots \\
(-1)^k \det (A^{(1, \ldots , k-1, j)}) \\
0 \\ 
\vdots \\
0 \\ 
(-1)^{k+1 }\det(A^{(1, \ldots , k)}) \\
0 \\ 
\vdots \\
0 
\end{pmatrix},$$
 where $ (-1)^{k+1}
\det(A^{(1, \ldots , k)})$ is the $j$-th entry.
\end{lem}

\begin{proof} The vectors $v_{k+1}, \ldots,v_n$ are obviously linearly independent,
they are $ n-k$ and we can easily see that they are in the kernel of $A$, so we conclude.
\end{proof}

\begin{cor} \label{bello} (1)
Let $K$ be a field and let $k,n \in \N\setminus \{0\}$. Let $A \in M(k \times n ,K)$ with $ n > k$ and $\rk(A)=k$.
For any  $j_1, \ldots, j_{k+1} $ in $\{1, \ldots, n\}$ with $j_1 < \ldots < j_{k+1} $, let $v_{j_1, \ldots, j_{k+1}}$ be the vector such that
\begin{itemize}
\item the $i$-th entry is equal to $0$ 
for every $i \neq j_1, \ldots, j_{k+1}$,
\item
 the  $j_l$-entry, for $l=1, \ldots, k+1$,
is equal to $$(-1)^l \det (A^{(j_1, \ldots, \hat{j_l}, \ldots, j_{k+1}) }).$$ 
 \end{itemize}
Then 
the kernel of $A$ is generated by the vectors 
 $v_{j_1, \ldots, j_{k+1}}$  for $j_1, \ldots, j_{k+1} $ elements of $\{1, \ldots, n\}$  with $j_1 < \ldots < j_{k+1} $.

  (2) 
Let $K$ be a field and let $m,n,k \in \N\setminus \{0\}$
with $n > k$. Let $A \in M(m \times n ,K)$
with $n > \rk(A) \geq k$.
For any 
 $i_1, \ldots,i_s $ in $\{1, \ldots, m\}$  with $i_1 < \ldots < i_s $, 
for any  $j_1, \ldots, j_{s+1} $ in $\{1, \ldots, n\}$  with $j_1 < \ldots < j_{s+1} $, let $v^{i_1, \ldots, i_s}_{j_1, \ldots, j_{s+1}}$ be the vector such that
\begin{itemize}
\item the $i$-th entry is equal to $0$ 
for every $i \neq j_1, \ldots, j_{s+1}$,
\item
 the  $j_l$-entry, for $l=1, \ldots, s+1$,
is equal to $$(-1)^l \det (A^{(j_1, \ldots, \hat{j_l}, \ldots, j_{s+1}) }_{(i_1, \ldots, i_s)}).$$ 
 \end{itemize}
Then the kernel of $A$ is generated by the vectors 
$v^{i_1, \ldots, i_s}_{j_1, \ldots, j_{s+1}}$
for  $s \in \{ k, \ldots,$  $\min\{m,n-1\} \}$, 
 $i_1, \ldots,i_s $ in $\{1, \ldots, m\}$  with $i_1 < \ldots < i_s $, 
$j_1, \ldots, j_{s+1} $ in $\{1, \ldots, n\}$  with $j_1 < \ldots < j_{s+1} $.
\end{cor}

\begin{thm} \label{Q-R2}
Let  $ p \geq q$ and let   
$\al =( [\underline{\alpha}_{i,j}, \overline{\alpha}_{i,j}])_{i,j}$ be a $p \times q$ interval matrix with 
$\underline{\alpha}_{i,j} \leq \overline{\alpha}_{i,j}$ and 
 $\underline{\alpha}_{i,j}, \overline{\alpha}_{i,j} \in \Q$ for any $i,j$. 
If there exists $ R \in \al$ with $\rk(R) =2$, then 
there exists $ Q \in \al \cap M(p \times q, \Q)$ with $\rk(Q) =2$.
\end{thm}

\begin{proof}
Let $a_i, b_i, c_j, d_j \in \R$ for $i=1,\ldots, 
p$ and  $j=1, \ldots, q$ such that 
$$r_{i,j} = a_i c_j + b_i d_j$$
for any $i,j$. Observe that we can easily suppose that, for any $j \in \{1, \ldots, q\}$, at least one of $ c_j $ and $d_j $ is nonzero (call this assumption ``assumption~($\ast$)'').

For any $i \in \{1, \ldots, p\}$   such that $\al_{i,j}$ is nondegenerate   for at least one   $j \in \{1, \ldots, q\}$  and   
for any $j \in \{1, \ldots, q \}$  such that  $\al_{i,j}$ is nondegenerate for at least one $i \in \{1, \ldots, p\}$,
 let $A_i, B_i, C_j, D_j$ be open neighbours respectively of $a_i, b_i, c_j, d_j$ such that 
$$ A_i C_j +  B_i D_j 	\subset \al_{i,j}$$
for any $(i,j) \in \{1, \ldots, p\} \times \{1, \ldots, q\}$ such that $\al_{i,j}$ is nondegenerate.  

Define  $$ \begin{array}{l}
T=\{ (i,j) \in \{1, \ldots, p\} \times \{1, \ldots, q\}    | \; \al_{i,j}\;  \mbox{\rm is degenerate}  \},  \vspace*{0.2cm}
\\ 
T_1=\{ i \in \{1, \ldots, p\}    | \; \exists  \; j \in \{1, \ldots, q\}  \;   \mbox{\rm s.t. } (i,j) \in T  \},
\vspace*{0.2cm}
\\
T_2=\{ j \in \{1, \ldots, q\}    | \; \exists  \; i  \in \{1, \ldots, p\}   \;   \mbox{\rm s.t. } (i,j) \in T  \} ,
\end{array}
 $$  
$$t_1 = |T_1| ,\hspace*{2cm} t_2 = |T_2|;$$

 for any $ i \in T_1 $, let $$T(i,\cdot)= \{ 
j \in \{1, \ldots, q\}    |  (i,j) \in T  \}$$
and, for 
any $ j \in T_2 $, let $$T(\cdot, j)= \{ 
i \in \{1, \ldots, p\}    |  (i,j) \in T  \}.$$

We can easily suppose that 
$T$ is nonempty, hence $t_1>0$ and $t_2>0$, and that
$$T_1= \{1, \ldots, t_1 \} ,
\hspace*{2cm} T_2= \{1, \ldots, t_2 \}.$$
Obviously, for any $(i,j) \in T$ 
\begin{equation} \label{ai}
a_i c_j + b_i d_j = \al_{i,j} ,
\end{equation}
where here $\al_{i,j}$  denotes one of the two (equal) endpoints of $\al_{i,j} $. So, if $(i,j) $ and $(i,h) $ are in $T$, we have:
$$ 
c_h ( \al_{i,j} -b_i \, d_j)=  
c_j ( \al_{i,h} -b_i \, d_h),$$ 
thus 
\begin{equation} \label{bi}
b_i ( c_j \,d_h - c_h \,d_j) = c_j \,\al_{i,h} -c_h \, \al_{i,j}.
\end{equation}
By (\ref{ai}), equation (\ref{bi})
 holds also if $(i,h ) \in T$ and  $c_h=0$, for
any $j =1,\ldots, q$.

From $(\ref{ai})$, we can deduce also 
that, if $(i,j)  $ and $(i,h) $ are in $T$, then
\begin{equation} \label{aibis}
a_i ( c_h \,d_j-c_j \,d_h ) = d_j \,\al_{i,h} -d_h \, \al_{i,j}.
\end{equation}


Moreover, by  (\ref{bi}),  if $(i,j) $, $(i,h) $ and $(i,k)$
are in $T$, 
then
$$ 
(c_j d_k -c_k d_j) (c_j \, \al_{i,h} - c_h \,\al_{i,j} )
= 
(c_j d_h -c_h d_j) (c_j \, \al_{i,k} - c_k \, \al_{i,j} ),
$$ that is
\begin{equation} \label{system}
\begin{array}{l}
 c_j (c_j \, \al_{i,h} -c_h \, \al_{i,j}) d_k +\\
-c_j (c_j \, \al_{i,k}   -c_k \, \al_{i,j}) d_h +\\
+[c_h (   c_j \, \al_{i,k}   -c_k \, \al_{i,j}  )
-c_k ( c_j \, \al_{i,h} -c_h \, \al_{i,j} )  ]
d_j =0 
\end{array}
\end{equation}
Let us consider the homogeneous linear system (S) in the unknowns $\delta_j$  for $j =1, \ldots,  t_2$
given by the equations  
\begin{equation} \label{systembis}
\begin{array}{l}
 \gamma_j (\gamma_j \, \al_{i,h} - \gamma_h \, \al_{i,j}) \delta_k +\\
-\gamma_j (\gamma_j \, \al_{i,k}   -\gamma_k \, \al_{i,j}) \delta_h +\\
+[ \gamma_h (   \gamma_j \, \al_{i,k}   -\gamma_k \, \al_{i,j}  ) -\gamma_k ( \gamma_j \, \al_{i,h} -\gamma_h \, \al_{i,j} ) ]
\delta_j =0 
\end{array}
\end{equation}
for any $i,j,k,h$ such that $(i,j) $, $(i,h) $ and $(i,k)$
are in $T$ and
the equations
\begin{equation} \label{systemtris}
\begin{array}{l} \al_{i,h}  \delta_k - \al_{i,k}   \delta_h =0 
\end{array}
\end{equation}
 for any $i,h,k$ such that 
$(i,h) $ and $(i,k)$
are in $T$ and $c_h=c_k=0$
(observe that the first equations are obtained from (\ref{system}) by replacing $c_i$ with $\gamma_i$ and $d_i$ with the unknown $\delta_i$).

Let us denote by $G_{(\al_{r,s})_{(r,s) \in T}, ( \gamma_j)_{j \in T_2}} $, or by $G_{(\al_{r,s}), ( \gamma_j)} $ for short,
the associated matrix, which has obviously $t_2$ columns.

If $\rk (G_{(\al_{r,s}), ( c_j)}) \geq 1 $, let  
$\overline{C}= \overline{C}_1 \times \ldots \times \overline{C}_{t_2}$  be 
a neighbourhood of 
$ (c_1,\ldots, c_{t_2})$ contained in $C_1 \times \ldots \times C_{t_2}$ 
  such that 
$\rk (G_{(\al_{r,s}), ( \gamma_j)}) \geq 1 $
for every $(\gamma_j)_j \in \overline{C}$.

 By Corollary \ref{bello},
 if   $1 \leq \rk(G_{(\al_{r,s}), ( \gamma_j)}) \leq t_2-1 $,
  we can see the kernel of  $G_{(\al_{r,s}), ( \gamma_j)} $  
for  $(\gamma_j)_j  \in \overline{C}$
as generated by some
vectors $$v_f  \left((\gamma_j)_{j \in T_2} \right) \;\; \; \mbox{\rm for } f=1, \ldots, g$$ (for some $g$) 
 whose entries are  polynomials in  $\al_{r,s}$  for $(r,s) \in T$ (which are fixed)
and $\gamma_j$ for $j \in T_2$.
By (\ref{system}) and (\ref{aibis}),  $(d_1, \ldots, d_{t_2}) $ satisfies 
both the equations (\ref{systembis}) and the equations (\ref{systemtris}) with $c_j$ instead of $\gamma_j$; moreover also 
$(c_1, \ldots, c_{t_2}) $ satisfies 
both the equations (\ref{systembis}) and the equations (\ref{systemtris}) with $c_j$ instead of $\gamma_j$; hence 
$\rk (G_{(\al_{r,s}), ( c_j)}) \leq t_2 - 1 $, because at least one of $(c_1, \ldots, c_{t_2}) $ and $(d_1, \ldots, d_{t_2}) $ must be nonzero by assumption ($\ast$);
moreover, since $(d_1, \ldots, d_{t_2}) $ satisfies 
both the equations (\ref{systembis}) and the equations (\ref{systemtris}) with $c_j$ instead of $\gamma_j$, we have that,  
if $\rk (G_{(\al_{r,s}), ( c_j)}) \geq 1 $, then
$$ \begin{pmatrix} 
d_1 \\
\vdots \\
d_{t_2}
\end{pmatrix}= \lambda_1
v_1  \left( (c_j)_{j \in T_2} \right) 
+ \ldots+ \lambda_g v_g \left((c_j)_{j \in T_2} \right) ,
$$
for some $\lambda_1, \ldots, \lambda_g \in \R$.

\underline{If $\rk (G_{(\al_{r,s}), ( c_j)}) \geq 1 $,}
choose 
\smallskip

(i) $\tilde{c}_j \in
\overline{C}_j  \cap \Q $ 
for any $j \in T_2$, 
\smallskip

(ii)
$(\tilde{\lambda}_1, \ldots, \tilde{\lambda}_g)$  in a neighbourhood 
of $(\lambda_1, \ldots, \lambda_g) $ and in  $\Q^g$,

in such a way that, if we
 define 
\begin{equation} \label{dtilde}
 \begin{pmatrix} 
\tilde{d}_1 \\
\vdots \\
\tilde{d}_{t_2}
\end{pmatrix}= \tilde{\lambda}_1
v_1  \left( (\tilde{c}_j)_{j \in T_2} \right) 
+ \ldots+ \tilde{\lambda}_g v_g \left( (\tilde{c}_j)_{j \in T_2} \right) ,
\end{equation}
we have that:

(a) $\tilde{c}_j=0$ if and only if $ c_j=0$;

(b) if $i, j,h$ are such that $(i,j), (i.h)  \in T$, $c_j, c_h, \al_{i,j}, \al_{i,h}$ are nonzero and 
$\det \begin{pmatrix} c_j & \al_{i,j} \\
 c_h & \al_{i,h} 
\end{pmatrix}=0$, then 
$\det \begin{pmatrix} \tilde{c}_j & \al_{i,j} \\
 \tilde{c}_h & \al_{i,h} 
\end{pmatrix}=0$;

(c)  if   
$ \det \begin{pmatrix} c_j &
 d_j  \\ c_h & d_h  
\end{pmatrix} \neq 0$, then 
$ \det \begin{pmatrix} \tilde{c}_j &
 \tilde{d}_j  \\ \tilde{c}_h & \tilde{d}_h  
\end{pmatrix} \neq 0$;

(d)  for any $ j\in T_2$, if  $d_j \neq 0 $, then 
$ \tilde{d}_j  \neq 0$;

(e) $ \tilde{d}_j \in D_j$  for any $ j\in T_2$.

Observe that, by  the choice of  $ (\tilde{c}_1, \ldots, \tilde{c}_{ t_2})$ we have done,
 the rank of  $G_{(\al_{r,s}), ( \tilde{c}_j)} $  is less than or equal to $t_2-1$:
 if $(\tilde{c}_1, \ldots, \tilde{c}_{t_2}) =0 $, then the equations (\ref{systembis}) with $\tilde{c}_j $ instead of $\gamma_j$ become trivial; moreover, by (a), we have that $(c_1, \ldots, c_{t_2}) =0 $, hence, by ($\ast$), 
  $(d_1, \ldots, d_{t_2}) \neq 0 $;
  by (\ref{aibis}), we have that  $(d_1, \ldots, d_{t_2}) $ satisfies equations (\ref{systemtris}) and then the system (S);
 if $(\tilde{c}_1, \ldots, \tilde{c}_{t_2})  \neq 0 $, the statement follows from the fact that 
   the transpose of $ (\tilde{c}_1,\ldots,
 \tilde{c}_{t_2})$ is in the  kernel of  $G_{(\al_{r,s}), ( \tilde{c}_j)} $.

Observe also that,
 obviously, $\tilde{d}_j \in \Q$ for any $j \in T_2$,
  since the $ \tilde{c}_j $ for $j \in T_2$ and the $ \tilde{\lambda}_f $ for $f=1,\ldots, g $ are in $\Q$.
\smallskip

\underline{If $\rk (G_{(\al_{r,s}), ( c_j)}) =0 $,}
choose 
\smallskip

(i) $\tilde{c}_j \in
C_j  \cap \Q $ 
for any $j \in T_2$, 
\smallskip

\smallskip
(ii)
$ \tilde{d}_j  \in D_j \cap \Q$  for any $ j\in T_2$,

in such way that

(a) $\tilde{c}_j=0$ if and only if $ c_j=0$;

(b) 
 if $i, j,h$ are such that $(i,j), (i,h)  \in T$,
   $c_j, c_h, \al_{i,j}, \al_{i,h}$ are nonzero and 
$\det \begin{pmatrix} c_j & \al_{i,j} \\
 c_h & \al_{i,h} 
\end{pmatrix}=0$, then 
$\det \begin{pmatrix} \tilde{c}_j & \al_{i,j} \\
 \tilde{c}_h & \al_{i,h} 
\end{pmatrix}=0$;

(c)  if 
$ \det \begin{pmatrix} c_j &
 d_j  \\ c_h & d_h  
\end{pmatrix} \neq 0$, then 
$ \det \begin{pmatrix} \tilde{c}_j &
 \tilde{d}_j  \\ \tilde{c}_h & \tilde{d}_h  
\end{pmatrix} \neq 0$;

(d)  for any $ j\in T_2$, if  $d_j \neq 0 $, then 
$ \tilde{d}_j  \neq 0$.

\medskip

{\bf Remark A.} {\em
Observe that, by our choices, if  
$\rk (G_{(\al_{r,s}), ( c_j)}) =0 $, then 
$\rk (G_{(\al_{r,s}), ( \tilde{c}_j)}) =0 $.}

This follows from the following remarks. 
\begin{itemize}
\item
Equations (\ref{systemtris}) do not depend on  $\gamma_j$.

\item Let us consider equations (\ref{systembis}); 
\subitem{-} if, for any $i,j,h,k$ 
with $j,h,k \in T(i, \cdot) $, we have that 
$c_j=c_h=c_k=0$, then,  by condition (a),
we have that   $\tilde{c}_j=\tilde{c}_h=\tilde{c}_k=0$, hence equations (\ref{systembis}) with $\tilde{c}_{l}$ instead
of $\gamma_l$ for every $l$ are trivial;
\subitem{-} if there exist $i,j,h,k$ 
with $j,h,k \in T(i, \cdot) $ and 
$c_j \neq 0$, then, since  $\rk (G_{(\al_{r,s}), ( c_j)}) =0 $, we must have 
$$\det \begin{pmatrix} c_j & \al_{i,j} \\
 c_h & \al_{i,h} 
\end{pmatrix}= \det \begin{pmatrix} c_j & \al_{i,j} \\
 c_k & \al_{i,k} 
\end{pmatrix}=0;$$
hence
\begin{equation} \label{fing}
\det \begin{pmatrix} \tilde{c}_j & \al_{i,j} \\
 \tilde{c}_h & \al_{i,h} 
\end{pmatrix}= \det \begin{pmatrix} \tilde{c}_j & \al_{i,j} \\
 \tilde{c}_k & \al_{i,k} 
\end{pmatrix}=0,
\end{equation}
 in fact: let us prove for example that 
$\det \begin{pmatrix} \tilde{c}_j & \al_{i,j} \\
 \tilde{c}_h & \al_{i,h} 
\end{pmatrix}= 0$:  

\subsubitem{{\Large $\cdot$}}
if $c_h=0$ then $\al_{i,h}$ must be zero (from $\det \begin{pmatrix} c_j & \al_{i,j} \\
 c_h & \al_{i,h} 
\end{pmatrix}=0$ and $c_j \neq 0$) and $\tilde{c}_h $ must be zero (by condition (a)), thus 
$\det \begin{pmatrix} \tilde{c}_j & \al_{i,j} \\
 \tilde{c}_h & \al_{i,h} 
\end{pmatrix}=0$;

\subsubitem{{\Large $\cdot$}}
if $c_h \neq 0$, then, from  $\det \begin{pmatrix} c_j & \al_{i,j} \\
 c_h & \al_{i,h} 
\end{pmatrix}=0$ and $c_j \neq 0$ we get that 
either $\al_{i,h}= \al_{i,j}=0$ or both    $\al_{i,h}$ and $ \al_{i,j}$ are nonzero;
if $\al_{i,h}= \al_{i,j}=0$, then obviously
$\det \begin{pmatrix} \tilde{c}_j & \al_{i,j} \\
 \tilde{c}_h & \al_{i,h} 
\end{pmatrix}=0$; if  both    $\al_{i,h}$ and $ \al_{i,j}$ are nonzero, then 
$\det \begin{pmatrix} \tilde{c}_j & \al_{i,j} \\
 \tilde{c}_h & \al_{i,h} 
\end{pmatrix}=0$ by condition (b).

\medskip

Moreover observe that (\ref{fing}) and the fact that $ \tilde{c}_j \neq 0$ imply that 
$\det \begin{pmatrix} \tilde{c}_k & \al_{i,k} \\
 \tilde{c}_h & \al_{i,h} 
\end{pmatrix}=0$.
Therefore, in every case, equations (\ref{systembis}) with $\tilde{c}_{l}$ instead
of $\gamma_l$ for every $l$ are trivial.

\end{itemize}

\smallskip

{\bf  Remark B.} {\em 
If,  for some $i,j,h$ with  $j,h \in T(i, \cdot) $ we have that   $ \det \begin{pmatrix} c_j &
 d_j  \\ c_h & d_h  
\end{pmatrix}=0$ (in particular, by (c), this holds if
$ \det \begin{pmatrix} \tilde{c}_j &
 \tilde{d}_j  \\ \tilde{c}_h & \tilde{d}_h  
\end{pmatrix}=0$), then 
$\det \begin{pmatrix} \tilde{c}_j & \al_{i,j} \\
 \tilde{c}_h & \al_{i,h} 
\end{pmatrix}=0$.}

 In fact:    
$ \det \begin{pmatrix} c_j &
 d_j  \\ c_h & d_h  
\end{pmatrix}=0$ implies, by (\ref{bi}), that 
$\det \begin{pmatrix} c_j & \al_{i,j} \\
 c_h & \al_{i,h} 
\end{pmatrix}=0$; thus 
$\det \begin{pmatrix} \tilde{c}_j & \al_{i,j} \\
 \tilde{c}_h & \al_{i,h} 
\end{pmatrix}=0$ (let $z$ be the cardinality of the nonzero  entries of 
 $\begin{pmatrix} c_j & \al_{i,j} \\
 c_h & \al_{i,h} 
\end{pmatrix}$; if $z=4$ our statement is true by (b), if $z \leq 2$, it is obviously  true; observe that  $z$ cannot be $3$).

\medskip
We have defined $\tilde{c}_j$ and $\tilde{d}_j$ for any $j \in T_2$. We want now to define 
$\tilde{a}_i$ and $\tilde{b}_i$ for any $i \in T_1$.
Let  $i \in T_1$. 
\begin{itemize}
\item If $|T(i,\cdot)| =1$, let $ T(i, \cdot)=\{ \overline{\jmath}(i) \}$; choose $ \tilde{a}_i $ and $\tilde{b}_i$ in $\Q$  such that 
\begin{equation} \label{Ti1}
 \tilde{a}_i  \tilde{c}_{\overline{\jmath}(
i)}   +\tilde{b}_i  \tilde{d}_{\overline{\jmath}(
i)} = \al_{i, \overline{\jmath}(i)}.
\end{equation}

\item Suppose  $|T(i,\cdot)| \geq 2$.

\subitem CASE 1:
if $  \tilde{c}_j \tilde{d}_h - \tilde{c}_h \tilde{d}_j $ is nonzero for some distinct $j ,h \in T(i,\cdot)$  (hence, in particular,  there exists $j \in T(i, \cdot)$ with  $\tilde{c}_j \neq 0$), define $\tilde{b}_i$ (by analogy with (\ref{bi}))  as follows:
\begin{equation} \label{tildebi}
\tilde{b}_i = \frac{\tilde{c}_j \,\al_{i,h} -\tilde{c}_h \, \al_{i,j}}{ \tilde{c}_j \tilde{d}_h - \tilde{c}_h \tilde{d}_j }
\end{equation}
(observe that the denominator  is nonzero by our assumption; moreover, it is a good definition, i.e. it does not depend on the choice of $j,h \in T(i , \cdot)$ such that $\tilde{c}_j \tilde{d}_h - \tilde{c}_h \tilde{d}_j \neq 0$, because the transpose of  $ (\tilde{d}_1, \ldots, 
\tilde{d}_{t_2})$ is  in $ Ker ( G_{(\al_{r,s}), ( \tilde{c}_j)}) $  (by (\ref{dtilde}) 
in case $\rk (G_{(\al_{r,s}), ( c_j)}) \geq 1 $
and  by Remark A in case $\rk (G_{(\al_{r,s}), ( c_j)}) =0 $), hence it satisfies the equations (\ref{systembis})  with $\tilde{c}_j$ instead of $ \gamma_j$);
  then define $\tilde{a}_i$ for $i=1, \ldots,
t_1$   by
\begin{equation} \label{tildeai}
\tilde{a}_i   = \frac{\al_{i,j} -  \tilde{b}_i  \tilde{d}_j}{ \tilde{c}_j}
\end{equation}
for any $j \in T(i, \cdot)$ with  $\tilde{c}_j \neq 0$;
it is a good definition by our definition 
of $ \tilde{b}_i$ and by Remark B, in fact: let $j, h \in T(i, \cdot) $ such that $\tilde{c}_j \neq 0 $ and $ \tilde{c}_h \neq 0$; we have to prove that  $  \frac{\al_{i,j} -  \tilde{b}_i  \tilde{d}_j}{ \tilde{c}_j} =  \frac{\al_{i,h} -  \tilde{b}_i  \tilde{d}_h}{ \tilde{c}_h}$; this is equivalent to $  
(\tilde{c}_j \tilde{d}_h - \tilde{c}_h \tilde{d}_j )
\tilde{b}_i =\tilde{c}_j \,\al_{i,h} -\tilde{c}_h \, \al_{i,j}$, which is true by the definition of  $ \tilde{b}_i$ in case $\tilde{c}_j \tilde{d}_h - \tilde{c}_h \tilde{d}_j \neq 0$ and by Remark B in case  $\tilde{c}_j \tilde{d}_h - \tilde{c}_h \tilde{d}_j = 0$.

\subitem CASE 2:
if $  \tilde{c}_j \tilde{d}_h - \tilde{c}_h \tilde{d}_j =0$  for every $j ,h \in T(i,\cdot)$, 
then, by Remark B, 
 we have that  
$\det \begin{pmatrix} \tilde{c}_j & \al_{i,j} \\
 \tilde{c}_h & \al_{i,h} 
\end{pmatrix}=0$  for every $j ,h \in T(i,\cdot)$.

\subsubitem
Case 2.1. If there exists $  j \in T(i, \cdot) $ such that 
$\tilde{c}_j \neq 0$, 
define $\tilde{b}_i$ to be any element of $B_i \cap \Q$ and 
$\tilde{a}_i $ as in $( \ref{tildeai})$ for any $j  
\in T(i, \cdot) $ such that 
$\tilde{c}_j \neq 0$  (it is well defined because 
$  \tilde{c}_j \tilde{d}_h - \tilde{c}_h \tilde{d}_j $ 
and $  \tilde{c}_j \al_{i,h} - \tilde{c}_h \al_{i,j} $  
are  zero hence $\tilde{b}_i( \tilde{c}_j \tilde{d}_h - \tilde{c}_h \tilde{d}_j)= 
  \tilde{c}_j \al_{i,h} - \tilde{c}_h \al_{i,j} $ for any $j ,h \in T(i, \cdot)$). 

\subsubitem
Case 2.2.
If $\tilde{c}_j = 0$ for any $  j \in T(i, \cdot) $,
then, by (a),  $c_j = 0$ for any $  j \in T(i, \cdot) $; hence  $d_j \neq 0$ for any $  j \in T(i, \cdot) $ by assumption~($\ast$); therefore, by condition (d), we have that
 $\tilde{d}_j \neq 0$ for any $  j \in T(i, \cdot) $; 
 define 
\begin{equation} \label{17}
 \tilde{b}_i= \frac{\al_{i,j}}{\tilde{d}_j}
 \end{equation}
 for any $j  \in T(i, \cdot) $  (it is a good definition
  because, as we have already said,  $ (\tilde{d}_1, \ldots ,
\tilde{d}_{t_2})$ satisfies the system (S)  with $\tilde{c}_j$ instead of $ \gamma_j$, in particular satisfies equations (\ref{systemtris})).
Moreover define $\tilde{a}_i$ to be any element of $A_i \cap \Q$.
\end{itemize}
 
By continuity, we can do the choices of the $
\tilde{c}_j$ and  of the $\tilde{\lambda}_f$
in case  $\rk (G_{(\al_{r,s}), ( c_j)}) \geq 1 $,
of the $
\tilde{c}_j$ and  of the $\tilde{d}_j$
in case  $\rk (G_{(\al_{r,s}), ( c_j)}) =0 $
in such way that:

- if $|T(i, \cdot)|=1$ we can choose the $\tilde{a}_i $ and the $\tilde{b}_i$ satisfying (\ref{Ti1}) respectively in $A_i $ and $B_i$,

- if $|T(i, \cdot)| \geq 2 $, the $\tilde{a}_i $ and the $\tilde{b}_i$ we have defined are  respectively in $A_i $ and $B_i$.

Finally, 
for any $ i \in \{1, \ldots, p\} \setminus T_1$, define $$ \tilde{a}_i = 
\left\{ \begin{array}{ll} a_i & \mbox {\rm if } a_i \in \Q 
\\ \mbox{\rm a point of } A_i \cap \Q & \mbox {\rm if } a_i \in \I,
 \end{array}   
\right. $$
$$\tilde{b}_i = 
\left\{ \begin{array}{ll} b_i & \mbox {\rm if } b_i \in \Q 
\\ \mbox{\rm a point of } B_i \cap \Q & \mbox {\rm if } b_i \in \I;
 \end{array}   
\right. 
$$  and,
for any $ j \in \{1, \ldots, q\} \setminus T_2$, define $$ \tilde{c}_j = 
\left\{ \begin{array}{ll} c_j & \mbox {\rm if } c_j \in \Q 
\\ \mbox{\rm a point of } C_j \cap \Q & \mbox {\rm if } c_j \in \I,
 \end{array}   
\right.$$ 
$$\tilde{d}_j = 
\left\{ \begin{array}{ll} d_j & \mbox {\rm if } d_j \in \Q 
\\ \mbox{\rm a point of } D_j \cap \Q & \mbox {\rm if } d_j \in \I.
 \end{array}   
\right. 
$$ 
We have 
$$
\tilde{a}_i  \tilde{c}_j + \tilde{b}_i  \tilde{d}_j  
\in \al_{i,j}$$  
for any $(i,j ) \not\in T $, since $\tilde{a}_i \in A_i$, 
 $\tilde{b}_i \in B_i$, $\tilde{c}_j \in C_j$, $\tilde{d}_j \in D_j$  for any 
$i=1, \ldots, p$ and $j=1,\ldots,q$. Moreover
$\tilde{a}_i$, 
 $\tilde{b}_i $, $\tilde{c}_j $, $\tilde{d}_j $ are in $\Q$  for any 
$i=1, \ldots, p$ and $j=1,\ldots,q$. 

Finally we prove that 
$$
\tilde{a}_i  \tilde{c}_j + \tilde{b}_i  \tilde{d}_j  
= \al_{i,j}$$
for any $(i,j ) \in T $.

In case $|T(i, \cdot)|=1$ the statement is true by (\ref{Ti1}).

Suppose $|T(i, \cdot)| \geq 2$.

 $\bullet$ 
If $\tilde{c}_j \neq 0$ we can be either in Case 1 or
in Case 2.1; in both cases 
the statement is true by our definition of $\tilde{a}_i$ (see (\ref{tildeai})).

 $\bullet$ 
Suppose  $\tilde{c}_j = 0$ and there exist $h,k \in T(i, \cdot)$ 
such that  
$ \tilde{c}_h \tilde{d}_k-  \tilde{c}_k \tilde{d}_h \neq 0$; hence there exists
$ h \in T(i, \cdot)$ with $ \tilde{c}_h \neq 0$ and we are in Case 1.
From the fact that    $\tilde{c}_j = 0$ we have, by (a), that $c_j=0$; therefore, by assumption~($\ast$), we have that $d_j \neq 0$; hence, by (d), we have that  $\tilde{d}_j \neq 0$; hence
$ \tilde{c}_j \tilde{d}_h-  \tilde{c}_h \tilde{d}_j= 
 -  \tilde{c}_h \tilde{d}_j \neq 0$ and
 the statement holds by our definition of $ \tilde{b}_i$ (see (\ref{tildebi})).

 
 $\bullet$ 
 Finally,
suppose   $\tilde{c}_j = 0$ and  $ \tilde{c}_h \tilde{d}_k-  \tilde{c}_k \tilde{d}_h = 0$ for any 
$h,k \in T(i, \cdot)$.  As before, this implies  $c_j=0$ (by (a)) and then $d_j \neq 0$ by assumption~($\ast$), and finally, by (d),  $\tilde{d}_j \neq 0$. From the fact that 
 $ \tilde{c}_h \tilde{d}_j- \tilde{ c}_j \tilde{d}_h = 0$ for any 
$h \in T(i, \cdot)$, we get that $\tilde{c}_h \tilde{d}_j =0 $ for any $h \in T(i, \cdot)$, hence 
$\tilde{c}_h=0 $ for any $h \in T(i, \cdot)$ and in this case (Case 2.2) we have defined $ \tilde{b}_i$ by  (\ref{17}).
Then  
 the statement is true
by (\ref{17}).

\medskip
Hence the $(p \times q)$-matrix $Q$  whose
entry $(i,j)$ is   $
\tilde{a}_i  \tilde{c}_j + \tilde{b}_i  \tilde{d}_j  $, for any 
$i=1, \ldots, p$ and $j=1,\ldots,q$, is a rational matrix of rank less than or equal to $2$ contained in $\al$. If $\rk(Q)=1$,
by changing an entry of $Q$ in an appropriate way,
we can get a rational matrix of rank $2$ contained in $\al$.
\end{proof}

\section{Rational realizations of the ranks $q-2$,
 $q-1$, $q$}

\begin{rem} \label{Q-Rq}
Let  $ p \geq q$ and let   
$\al =( [\underline{\alpha}_{i,j}, \overline{\alpha}_{i,j}])_{i,j}$ be a $p \times q$ interval matrix with 
$\underline{\alpha}_{i,j} \leq \overline{\alpha}_{i,j}$ and 
 $\underline{\alpha}_{i,j}, \overline{\alpha}_{i,j} \in \Q$ for any $i,j$. 
Suppose there exists $ A \in \al$ with $\rk(A) =q$; then obviously
there exists $ B \in \al \cap M(p \times q, \Q)$ with $\rk(B) =q$.
\end{rem}

\begin{nota}
In the proof of the following theorem we will use the following notation for any $b \in \R$: $$b \wr v = \left\{
\begin{array}{lll}  b & \mbox{\rm if} & b \in \Q, \\
v & \mbox{\rm if} & b \in \I , \end{array} \right.$$ 
where $v$ can be a real number or an interval.     

Moreover, if we have two systems of equations,  $(P)$ and $(M)$, in the same unknowns, 
we call $(PM)$ the system  ``union'' of the 
two systems, i.e. the system given both by the equations of $(P)$ and the equations of $(M)$.
\end{nota}

\begin{rem} \label{sistemi}
(i) If a linear system with rational entries has a 
solution $c$ and $V$ is a neighbourhood of $c$, then there is a solution of the system contained in $V$ and with rational entries.

(ii) Let  $(S_t)$ be a linear system whose entries  depend linearly on a parameter $t\in \R^n$. Let $ \overline{t} \in \R^n$. If  the system $S_{\overline{t}}$  has a solution $b$, then for every neighbourhood $U$ of $b$, there exists a neighbourhood $V$ of $\overline{t}$ such that if $t \in V$, $S_t$ is solvable  and the dimension of the solution space of $S_{t} $  is equal to  the dimension of the solution space of $S_{\overline{t}} $, then
there is a solution of $S_t$ in $U$.
\end{rem}

\begin{thm} \label{Q-R111}
Let  $ p \geq q$ and let   
$\al =( [\underline{\alpha}_{i,j}, \overline{\alpha}_{i,j}])_{i,j}$ be a $p \times q$ interval matrix with 
$\underline{\alpha}_{i,j} \leq \overline{\alpha}_{i,j}$ and 
 $\underline{\alpha}_{i,j}, \overline{\alpha}_{i,j} \in \Q$ for any $i,j$. 
If there exists $ A \in \al$ with $\rk(A) \leq q-2$, then 
there exists $ B \in \al \cap M(p \times q, \Q)$ with $\rk(B) \leq q-2$.
\end{thm}

\begin{proof}
We can suppose that $A^{(q-1)}, A^{(q)} \in \langle A^{(1)}, \ldots  A^{(q-2)} \rangle $; 
let   
\begin{equation} \label{iniziale1}
A^{(q-1)}  = b_1 A^{(1)}+ \ldots+ b_{q-2}  A^{(q-2)}  
\end{equation}
for some $b_1,\ldots , b_{q-2} \in \R $ and let
\begin{equation} \label{iniziale2}
A^{(q)}  = c_1 A^{(1)}+ \ldots+ c_{q-2}  A^{(q-2)}  
\end{equation} 
for some $c_1,\ldots , c_{q-2} \in \R $. We can suppose that, for any $j=1, \ldots , q-2$, 
 either $b_j \neq 0 $ or $c_j \neq 0$ (call this assumption ($\ast \ast$)).

Up to swapping rows and columns, we can also suppose that 
$\al_{i,q-1}, \al_{i,q} $ are nondegenerate  for $i=1,  \ldots k$, while, for $i=k+1, \ldots p$, at least one of $\al_{i,q-1}$ and $\al_{i,q} $  is degenerate.

For any $ i=1,\ldots ,p$, we define:
$$N_i =\{ j \in \{1,\ldots , q-2 \} | \; a_{i,j} \in \I \}.$$

Moreover,  let us define $$ B=\{j \in \{1,\ldots, q-2\}| \; \; b_j \in \I  \},$$ 
  $$
   C=\{j \in \{1,\ldots, q-2\}| \; \; c_j \in \I  \}.$$

Finally we can  suppose that $a_{i,j} \in \Q $
for any $i=1,\ldots , k$ and $j=1,\ldots, q-2$; in
 fact: if 
 for some $i \in \{1, \ldots , k \}$ the set $N_i$ is nonempty, we have that for any $j \in N_i $ the entry $\al_{i,j}$ is nondegenerate
 (since has rational endpoints and contains $a_{i,j}$, which is irrational); so there exist  neighbourhoods $U_{i,j}$ of $a_{i,j}$ 
 contained in $\al_{i,j}$  for  any $j \in N_i $ 
 such that 
$$ \sum_{j \in \{1, \ldots , q-2\}} b_j  \cdot (a_{i,j} \wr U_{i,j}) \subset  \al _{i,q-1}  $$
and 
$$ \sum_{j \in \{1, \ldots , q-2\} } c_j  \cdot (a_{i,j} \wr U_{i,j}) \subset  \al _{i,q}  ;$$
hence, for any $j \in N_i$ we can change the entry $a_{i,j} $ into an element $ \tilde{a}_{i,j}$  of $ U_{i,j} \cap \Q$ and the 
entries $a_{i,q-1}$ and  $a_{i,q}$ respectively into
$$ \sum_{j \in \{1, \ldots , q-2\}} b_j \cdot (a_{i,j} \wr \tilde{a}_{i,j}), \hspace*{1cm}\mbox{\rm and} \hspace*{1cm}\sum_{j \in \{1, \ldots , q-2\} } c_j \cdot (a_{i,j} \wr \tilde{a}_{i,j}); $$
   in this way 
we get  again a matrix with the last two columns 
in the span of the first $q-2$ columns and such that 
the first $q-2$ entries of each of the first $k$ rows of  the matrix are  rational.  So  we can  suppose that $a_{i,j} \in \Q $
for any $i=1,\ldots , k$ and $j=1,\ldots, q-2$.

Moreover define  $$
\begin{array}{r}X^b=\{i \in \{k+1 , \ldots, p\} | \; 
\al_{i,q-1} \; \mbox{\rm is degenerate, }
\al_{i,q}  \; \mbox{\rm is nondegenerate   } \\
\mbox{\rm and either }
 N_i  = \emptyset \;
\mbox{\rm or } b_j =0 \;\forall j \in N_i
\}.
\end{array}
$$
$$
\begin{array}{r}X^c=\{i \in \{k+1 , \ldots, p\} | \; 
\al_{i,q-1} \; \mbox{\rm is nondegenerate, }
\al_{i,q}  \; \mbox{\rm is degenerate   } \\
\mbox{\rm and either }
 N_i  = \emptyset \;
\mbox{\rm or } c_j =0 \;\forall j \in N_i
\}.
\end{array}
$$
$$
\begin{array}{r}Y^b=\{i \in \{k+1 , \ldots, p\} | \; 
\al_{i,q-1} \; \mbox{\rm is degenerate, }
\al_{i,q}  \; \mbox{\rm is degenerate   } \\
\mbox{\rm and either }
 N_i  = \emptyset \;
\mbox{\rm or } b_j =0 \;\forall j \in N_i
\}.
\end{array}
$$
$$
\begin{array}{r}Y^c=\{i \in \{k+1 , \ldots, p\} | \; 
\al_{i,q-1} \; \mbox{\rm is degenerate, }
\al_{i,q}  \; \mbox{\rm is degenerate   } \\
\mbox{\rm and either }
 N_i  = \emptyset \;
\mbox{\rm or } c_j =0 \;\forall j \in N_i
\}.
\end{array}
$$
\medskip
We want now to define some neighbours $Z^i_j$ of $b_j$, $V^i_j$ of $c_j$ and $U_{i,j}$ of $a_{i,j}$ for some $j \in \{1, \ldots, q-2\}$ 
and $i \in \{1, \ldots, p\}$.
\medskip

{\bf CASE 0}: Let $i \in \{1 , \ldots, k\} $. 
Hence $\al_{i,q-1}$ and $\al_{i,q}$ are nondegenerate.

Choose 

 -  for any $j  \in B$,
a  neighbourhood $Z_j^i$ of $b_j$,

 - for any $j \in C$,         
  a neighbourhood $V_j^i$ of $c_j$, 
  
such that the following two conditions hold:
$$ \sum_{j \in \{1,\ldots, q-2\}}  (b_j \wr Z^i_j) \cdot a_{i,j} 
\subset \al_{i,q-1},$$ 
$$ \sum_{j \in \{1,\ldots, q-2\}} (c_j \wr V^i_j) \cdot a_{i,j} 
\subset \al_{i,q}.$$

\medskip

{\bf CASE 1}:  Let $i \in \{k+1 , \ldots, p\} $ such that 
$\al_{i,q-1}$ is degenerate and 
 $\al_{i,q}$ is nondegenerate. 
 
 $\bullet $ SUBCASE 1.1:
 $i \not \in  X^b$. 
 
 Hence $N_i \neq \emptyset $ and 
 there exists $\overline{\jmath}(i) \in N_i $ such that $b_{\overline{\jmath}(i)} \neq 0 $.

Choose

 -  for any $j  \in B$,
a  neighbourhood $Z_j^i$ of $b_j$,

 - for any $j \in C$,         
  a neighbourhood $V_j^i$ of $c_j$, 
  
 -  for any 
$j  \in N_i  $, a neighbourhood $U_{i,j}$ of $a_{i,j}$
contained in $\al_{i,j}$,

such that the following three conditions  hold: 

if $ b_{\overline{\jmath}(i)} \in \I$, we have that 
 $Z_{\overline{\jmath }(i)}^i$  is  contained either  in $\R_{<0}$ or in  $\R_{>0}$;

 \begin{equation} \label{tar}  \sum_{j \in \{1, \ldots, q-2\}} (c_j \wr V^i_j) \cdot  (a_{i,j} \wr U_{i,j}) 
  \subset    \al _{i,q};  \end{equation}

 \begin{equation} \label{tar2}
 -\frac{1}{(b_{\overline{\jmath}(i)} \wr Z^i_{\overline{\jmath}(i)})}
 \left[ 
\sum_{
j \in \{1, \ldots, q-2\}  \setminus 
\{\overline{\jmath}(i)\}
}
( b_j \wr
Z^i_j )\cdot (a_{i,j}\wr U_{i,j})  -a_{i,q-1} 
\right] 
 \subset 
U_{i, \overline{\jmath}(i)}.
\end{equation}

\medskip

 $\bullet $  SUBCASE 1.2: $i  \in  X^b$ and  $N_i \neq \emptyset$.

Choose

 - for any $j \in C$,         
  a neighbourhood $V_j^i$ of $c_j$,

  - for any 
$j  \in N_i  $, a neighbourhood $U_{i,j}$ of $a_{i,j}$
contained in $\al_{i,j}$,

such that  (\ref{tar})  holds.

 $\bullet $  SUBCASE 1.3: $i  \in  X^b$ and  $N_i = \emptyset$.

Choose

 - for any $j \in C$,         
  a neighbourhood $V_j^i$ of $c_j$,

such that  (\ref{tar})  holds with $a_{i,j}$ instead of $ (a_{i,j} \, \wr \, U_{i,j})$ (observe that $a_{i,j} \in \Q$ for any $j \in \{1,\ldots , q-2\}$).

\medskip

{\bf CASE 2}: 
Let $i \in \{k+1 , \ldots, p\} $ such that 
$\al_{i,q-1}$ is nondegenerate and 
 $\al_{i,q}$ is degenerate.

Analogous to Case 1 by swapping $q-1$ with $q$ and $b$ with $c$.

\medskip

{\bf CASE 3}:   Let $i  \in \{k+1, \ldots, p\}$ be such that $\al_{i,q-1}$ and 
 $\al_{i,q}$ are degenerate. 

\medskip
SUBCASE 3.1: $ i \not \in Y^b \cup Y^c$ and
there does not exist $j
 \in   N_i $ such that $b_j \neq 0 $ and $c_j \neq 0$;
 by the assumption ($\ast \ast$),
  there exist $ \overline{\jmath }(i)$ and  $\hat{\jmath}(i)$ 
 in  $ N_i $ such that $$b_{  \overline{\jmath }(i) } \neq 0, \;\;\;\; 
c_{  \overline{\jmath }(i) } = 0 ,  \;\;\;\; \;\;\;
c_{\hat{\jmath}(i)} \neq 0, \;\;\;\; 
b_{\hat{\jmath}(i)} = 0.$$

  We consider:

   for any $j  \in B$,
  neighbourhoods $Z_j^i$ of $b_j$,

  for any $j \in C$,         
  neighbourhoods $V_j^i$ of $c_j$, 
  
   for any 
$j  \in N_i  \setminus  \{\overline{\jmath}(i), \hat{\jmath}(i)
\}$, neighbourhoods $U_{i,j}$ of $a_{i,j}$
contained in $\al_{i,j}$

such that the following four conditions hold: 

if $ b_{\overline{\jmath}(i)} \in \I$, we have that 
 $Z_{\overline{\jmath }(i)}^i$  is  contained either  in $\R_{<0}$ or in  $\R_{>0}$;

if $ c_{\hat{\jmath}(i)} \in \I$, we have that 
 $V_{\hat{\jmath}(i)}^i$  is  contained either  in $\R_{<0}$ or in  $\R_{>0}$;

\begin{equation} \label{ar4}
  -\frac{1}{ b_{\overline{\jmath}(i)} \wr Z^i_{\overline{\jmath}(i)}}
 \left[ 
\sum_{j \in  \{1, \ldots,  q-2\} \setminus \{ \hat{\jmath}(i)
\overline{\jmath}(i)
\}} (b_j \wr
Z^i_j) \cdot ( a_{i,j} \wr U_{i,j})   -a_{i,q-1} 
\right] 
 \subset 
\al _{i, \overline{\jmath}(i)}; 
\end{equation}

\begin{equation} \label{ar1}
 -\frac{1}{c_{\hat{\jmath}(i)} \wr V^i_{\hat{\jmath}(i)}}
 \left[ 
\sum_{j \in  \{1, \ldots,  q-2\} \setminus \{ \hat{\jmath}(i)
\overline{\jmath}(i)
\}} (c_j \wr
V^i_j) \cdot ( a_{i,j} \wr U_{i,j}) 
  -a_{i,q} 
\right] 
 \subset 
\al _{i, \hat{\jmath}(i)}; 
\end{equation}

\medskip
SUBCASE 3.2: $ i \not \in Y^b \cup Y^c$ and
there exists $ \overline{\jmath }(i) \in 
  N_i$ such that $b_{\overline{\jmath }(i)} \neq 0$ and 
  $c_{\overline{\jmath }(i)} \neq 0$.

  We consider:

   for any $j  \in B$,
  neighbourhoods $Z_j^i$ of $b_j$, 
  
  for any $j \in C$,         
  neighbourhoods $V_j^i$ of $c_j$, 
  
   for any 
$j  \in N_i  \setminus  \{\overline{\jmath}(i)\}$, neighbourhoods $U_{i,j}$ of $a_{i,j}$
contained in $\al_{i,j}$

such that the following four conditions hold:

if $ b_{\overline{\jmath}(i)} \in \I$, we have that 
 $Z_{\overline{\jmath }(i)}^i$  is  contained either  in $\R_{<0}$ or in  $\R_{>0}$;

if $ c_{\overline{\jmath}(i)} \in \I$, we have that 
 $V_{\overline{\jmath }(i)}^i$  is  contained either  in $\R_{<0}$ or in  $\R_{>0}$;

\begin{equation} \label{star3}
  -\frac{1}{ b_{\overline{\jmath}(i)} \wr Z^i_{\overline{\jmath}(i)}}
 \left[ 
\sum_{j \in  \{1, \ldots,  q-2\} \setminus \{
\overline{\jmath}(i)
\}} (b_j \wr
Z^i_j) \cdot ( a_{i,j} \wr U_{i,j})   -a_{i,q-1} 
\right] 
 \subset 
\al _{i, \overline{\jmath}(i)}; 
\end{equation}

\begin{equation} \label{star4}
  -\frac{1}{ c_{\overline{\jmath}(i)} \wr V^i_{\overline{\jmath}(i)}}
 \left[ 
\sum_{j \in  \{1, \ldots,  q-2\} \setminus \{
\overline{\jmath}(i)
\}} (c_j \wr
V^i_j) \cdot ( a_{i,j} \wr U_{i,j})   -a_{i,q} 
\right] 
 \subset 
\al _{i, \overline{\jmath}(i)}; 
\end{equation}

\medskip
SUBCASE 3.3: $ i  \in Y^c \setminus Y^b$.

In this case we must have: $N_i \neq \emptyset$, $c_j=0 $ $\forall j \in N_i$ and there exists $ \overline{\jmath}(i)
\in N_i$ such that $b_{\overline{\jmath}(i)} \neq 0$.

  We consider:

   for any $j  \in B$,
  neighbourhoods $Z_j^i$ of $b_j$, 
  
   for any 
$j  \in N_i  \setminus  \{\overline{\jmath}(i)\}$, neighbourhoods $U_{i,j}$ of $a_{i,j}$
contained in $\al_{i,j}$

such that:

if $ b_{\overline{\jmath}(i)} \in \I$, we have that 
 $Z_{\overline{\jmath }(i)}^i$  is  contained either  in $\R_{<0}$ or in  $\R_{>0}$ and
(\ref{star3}) holds.

\medskip
SUBCASE 3.4: $ i  \in Y^b \setminus Y^c$. 

Analogous to the previous subcase. 

\medskip

Finally observe that if $i \in Y^b \cap Y^c$ we must have $N_i = \emptyset$ since, by  assumption ($\ast \ast$), there does not exist $j$ such that $b_j=c_j=0$.
In this case we do not give at the moment any definition.

\medskip {\em Definiition.}
$\bullet $ For any $j \in B$, let $ \beta(j) $ be the set of the $i
\in \{1, \ldots, p\}$ such that we have chosen
$Z^i_j$.
   
   $\bullet$
For any $j \in C$, let $ \gamma(j) $ be the set of the $i 
\in \{1, \ldots, p\}$ such that we have chosen $V^i_j$.  

\medskip

  \underline{{\em Choice of the $\tilde{b}_j$ for
  any $j \in B $ 
  and of the $\tilde{c}_j$ for
  any $j \in C $ in case $B \cup C \neq \emptyset$.   
  }}
  
  If $X^b  \cup Y^b = \emptyset$, then, for any $j \in B$, choose $ \tilde{b}_j $ in the set $$\left( \cap_{i \in \beta(j) } Z^i_j  \right) \cap \Q$$
  (observe that in this case we have $ \beta(j) \neq \emptyset$).

  If $X^c  \cup Y^c = \emptyset$, then, for any $j\in C$, choose $ \tilde{c}_j $ in the set $$\left( \cap_{i \in \gamma(j) } V^i_j  \right) \cap \Q$$
  (observe that in this case we have $ \gamma(j) \neq \emptyset$).
    
    If $X^b \cup Y^b  \neq \emptyset$,
define $(M)$ to be       the linear system  given by the equations 
\begin{equation} \label{sei}
\sum_{j \in \{1,\ldots , q-2\}} (b_j \wr
\tilde{b}_j)\cdot a_{i,j} = a_{i,q-1},
\end{equation}
  for $i \in X^b \cup Y^b $,
 in the unknowns $\tilde{b}_j$ for 
  $j\in B$.
Observe that the linear system has rational coefficients and it is certainly solvable since $b:=(b_j)_{j \in  B}$ is a  solution.  
   


    If $X^c \cup Y^c \neq \emptyset$, define $(N)$ to be the linear system  given by the equations 
\begin{equation} \label{sei}
\sum_{j \in \{1,\ldots , q-2\}} (c_j \wr
\tilde{c}_j)\cdot a_{i,j} = a_{i,q},
\end{equation}
  for $i \in X^c \cup Y^c $,
 in the unknowns $\tilde{c}_j$ for 
  $j\in C$.
Observe that the linear system has rational coefficients and it is certainly solvable since $c:=(c_j)_{j \in  C}$ is a  solution.  

\smallskip

Moreover,  define
$W$ to be the set
 $$\begin{array}{r}
 \{i \in \{1, \ldots, p\}| \;
\al_{i,q-1}, \al_{i,q} \; \mbox{\rm are degenerate, } i \not \in Y^b \cup Y^c, \vspace*{0.2cm} \\  
\hspace*{4cm}
\exists   \overline{\jmath }(i) \in 
  N_i \; \mbox{\rm s.t. } b_{\overline{\jmath }(i)} 
  c_{\overline{\jmath }(i)} \neq 0\}
  \end{array}$$
(that is $W$ is the set of the $i$ in Case 3.2)
and,
if $ B \cup C
\neq \emptyset$,  consider, for $ i \in W$,
the following  equation:
\begin{equation} \label{setot}
\begin{array}{c} 
-\frac{1}{b_{\overline{\jmath}(i)} \wr \tilde{b}_{\overline{\jmath}(i)}} \left[ 
\sum_{j \in  \{1, \ldots, q-2\} \setminus 
\{  \overline{\jmath}(i) \}
}( b_j \wr 
\tilde{b}_j)
\cdot ( a_{i,j} \wr \tilde{a}_{i,j})  -a_{i,q-1} 
\right] =\\
 -\frac{1}{c_{\overline{\jmath}(i)} \wr \tilde{c}_{\overline{\jmath}(i)}} \left[ 
\sum_{j \in  \{1, \ldots, q-2\} \setminus 
\{ \overline{\jmath}(i) \}
}( c_j \wr 
\tilde{c}_j)
\cdot ( a_{i,j} \wr \tilde{a}_{i,j})  -a_{i,q} 
\right], 
\end{array}
\end{equation}
 where  the unknowns are the 
$ \tilde{a}_{i,j}$ for
  $j \in N_i \setminus \{\overline{\jmath}(i)\}$.
 It is obviously equivalent to  
  
\begin{equation} \label{setot2}
\begin{array}{c} 
\sum_{j \in N_i \setminus \{\overline{\jmath}(i)\}}
\left(
\frac{b_j \wr \tilde{b}_j }{b_{\overline{\jmath}(i)} \wr \tilde{b}_{\overline{\jmath}(i)}} 
- 
\frac{c_j \wr \tilde{c}_j }{c_{\overline{\jmath}(i)} \wr \tilde{c}_{\overline{\jmath}(i)}} 
\right) \tilde{a}_{i,j}= 
\\ 
-\sum_{j \not\in N_i }
\left(
\frac{b_j \wr \tilde{b}_j }{b_{\overline{\jmath}(i)} \wr \tilde{b}_{\overline{\jmath}(i)}} 
- 
\frac{c_j \wr \tilde{c}_j }{c_{\overline{\jmath}(i)} \wr \tilde{c}_{\overline{\jmath}(i)}} 
\right) a_{i,j}+
\frac{a_{i,q-1}}{b_{\overline{\jmath}(i)} \wr \tilde{b}_{\overline{\jmath}(i)}}
- 
\frac{a_{i,q}}{c_{\overline{\jmath}(i)} \wr \tilde{c}_{\overline{\jmath}(i)}}
\end{array}
\end{equation}
  
 We would like to find $b:=(\tilde{b}_j)_{ j \in B}$ and $c:= (\tilde{c}_j)_{j \in C}$ 
such that they satisfy $(M)$ and $(N)$,  
   the equation (\ref{setot2}) in the unknowns  
$ \tilde{a}_{i,j}$ for
  $j \in N_i \setminus \{\overline{\jmath}(i)\}$ is solvable for any $i \in W$ and, finally, we can find a solution of  (\ref{setot2})  ``arbitrarily near'' to  
$ (a_{i,j})_{j \in N_i \setminus \{\overline{\jmath}(i)\}}$.
  
Let   $(\ref{setot2})_{b,c}$ be the equation we get from (\ref{setot2}) by replacing  $\tilde{b}_j$ with $b_j$ for every $j \in B$   and  $\tilde{c}_j$ with $c_j$ for every $j \in C$; observe that $(\ref{setot2})_{b,c}$  is true, because, 
   by replacing in (\ref{setot2}) $\tilde{b}_j$ with $b_j$ for every $j \in B$   and  $\tilde{c}_j$ with $c_j$ for every $j \in C$, both members of  (\ref{setot2}) become equal to $a_{i, \overline{\jmath}(i)}$.
  
Let $G$ be the set of the $i \in W$ such that 
at least one of the coefficients of $\tilde{a}_{i,j}$ for $ j \in N_i \setminus \{\overline{\jmath} (i)\}  $ in $(\ref{setot2})_{b,c}$ is nonzero.
  Let $Z$ and $V$ be neighbours respectively of $b$ and $c$ such that 
such coefficients  are nonzero for any $(\tilde{b}, \tilde{c}) \in Z \times V$ and any $i \in G$.
  
Let us call $(P)$   the linear system in the unknowns $\tilde{b}_j$ for $j \in B$  with entries depending on $\tilde{c}_j$  for $j \in C$ given by imposing  
the second member of $(\ref{setot2})$ 
 equal to $0$ for any $i \in W \setminus G$, i.e. given by the equations
  $$ (b_{\overline{\jmath}} \wr  \tilde{b}_{\overline{\jmath}})
  \left(a_{i,q} - \sum_{j \not\in  N_i }   (c_j \wr  \tilde{c}_j) a_{i,j} \right) + \sum_{j \not\in  N_i } (  b_j \wr  \tilde{b}_j)
(  c_{\overline{\jmath}} \wr  \tilde{c}_{\overline{\jmath}})
a_{i,j} -( c_{\overline{\jmath}} \wr  \tilde{c}_{\overline{\jmath}})  a_{i,q-1}=0$$
for any $i \in W \setminus G$, and given 
by the  equations
$$ 
\frac{b_j \wr \tilde{b}_j }{b_{\overline{\jmath}(i)} \wr \tilde{b}_{\overline{\jmath}(i)}} 
- 
\frac{c_j \wr \tilde{c}_j }{c_{\overline{\jmath}(i)} \wr \tilde{c}_{\overline{\jmath}(i)}} 
 =0$$ 
 for  $i \in W \setminus G$ and $j \in N_i \setminus 
\{\overline{\jmath}(i)\}$.

Let $(P_{c})$ be the system we get from $(P)$ by replacing $\tilde{c}$ with $ c$.

Let us call 
$(X)$ the linear system in the unknowns $\tilde{c}_j$ given by the equality :

\smallskip

\hspace*{2cm}
``rank of the incomplete matrix of $(PM)$=

\hspace*{2cm}
= rank of the incomplete matrix of $(P_c M)$=

\hspace*{2cm} = rank of the complete matrix of $(PM)$''.
  
  \smallskip
  
REMARK A. {\em
Observe that $c$ is a solution of $(X)$,   
in fact, $b$ is a solution   $(P_c, M)$. Moreover $c$ is a solution of $(N)$.
So $c$ is a solution of $(XN)$.}

\smallskip

REMARK B. {\em 
Observe also that, if $ b_{\overline{\jmath}}\in \I$,
 then all the columns of the complete matrix associated to the linear system $(P)$, apart from the column corresponding to $\tilde{b}_{\overline{\jmath}}$ 
 are multiple of 
$ c_{\overline{\jmath}} \wr \tilde{c}_{\overline{\jmath}}$ 
and not depending on the other 
$ c_{j} \wr \tilde{c}_{j}$'s  
and  it is easy to see  that the system $(X)$ is linear in 
 the $\tilde{c}_j$. 
 
 Also if $ b_{\overline{\jmath}}\in \Q$, we can conclude easily that the system $(X)$ is linear in 
 the $\tilde{c}_j$.}

\medskip

Suppose $C\neq \emptyset $.
By Remarks A and B,
the system  $(XN)$ is a linear system in the unknowns $\tilde{c}_j$ with rational entries and  $c$ 
is in its solution set. Hence, by (i) of Remark \ref{sistemi}, we can find a rational solution $\hat{c}:=(\hat{c}_j )_{j \in C}$ with $\hat{c}_j 
\in \cap_{i \in \gamma(j)} V^i_j $ for every $j$,
$\hat{c}_{\overline{\jmath}(i)} \neq 0$  for every $i \in W$, 
 and $\hat{c} \in V$ if $G \neq \emptyset$.

Let $(P_{\hat{c}})$ be the system we get from $(P)$ by replacing $\tilde{c}$ with $ \hat{c}$.
The linear system $(P_{\hat{c}}, M)$ 
  in the unknowns $\tilde{b}_j$
  for $j \in B$ is solvable 
  by our choice of $ \hat{c}$ 
  and the dimension of its solution set is equal to the dimension of the solution set of 
   $(P_{c}, M)$ 
  (in fact, $\hat{c}$ is a solution of  $(X)$); moreover it has rational entries, so it has a rational solution $ \hat{b} :=(\hat{b}_j )_{j \in B}$ by (i) of Remark \ref{sistemi}. 
  Moreover, by (ii) of Remark \ref{sistemi},
  we can choose $\hat{c}$ so that $ \hat{b}
  \in \cap_{i \in \beta(j)} Z^i_j $, $\hat{b}_{\overline{\jmath}(i)} \neq 0$ for every $i \in W$,  and $ \hat{b} \in Z$ if $ G \neq \emptyset$.
  
  If $C = \emptyset $, take $\hat{c}=c$ and argue analogously.
  
  The couple $(\hat{b}, \hat{c}) $ satisfies 
 $ (P)$ (precisely $\hat{b}$ satisfies $(P_{\hat{c}})$). So, if we replace
$(\tilde{b}, \tilde{c})$ with  $(\hat{b}, \hat{c}) $    in (\ref{setot}), we get a solvable  equation in the $\tilde{a}_{i,j}$.  Moreover 
$\hat{b}$ satisfies $(M)$  and 
$\hat{c}$ satisfies $(N)$. We choose $(\hat{b}, \hat{c}) $
for
$(\tilde{b}, \tilde{c})$.

If $C \neq \emptyset$ and $B = \emptyset$, we argue analogously.

\medskip
\underline{{\em Choice of the $\tilde{a}_{i,j} 
$ for
 $i \in \{k+1, \ldots , p\}$ and   $j \in N_i  $.}}

CASE 1, i.e.
$\al_{i,q-1}$ is degenerate and 
 $\al_{i,q}$ is nondegenerate.
\smallskip

  SUBCASE 1.1:  
 $i \not \in  X^b$
 (hence $N_i \neq \emptyset $ and 
 there exists $\overline{\jmath}(i) \in N_i $ such that $b_{\overline{\jmath}(i)} \neq 0 $).

Choose 
$\tilde{a}_{i,j} \in U_{i,j} \cap \Q$ for any $j \in N_i \setminus \{\overline{\jmath}(i)\}$ and define $ \tilde{a}_{i, \overline{\jmath}(i) }$ to be 
\begin{equation} \label{sette}
 -\frac{1}{b_{\overline{\jmath}(i)} \wr \tilde{b}_{\overline{\jmath}(i)}} \left[ 
\sum_{j \in  \{1, \ldots, q-2\} \setminus 
\{ \overline{\jmath}(i) \}
}( b_j \wr 
\tilde{b}_j)
\cdot ( a_{i,j} \wr \tilde{a}_{i,j})  -a_{i,q-1} 
\right]. 
\end{equation}

 By $( \ref{tar2})$, we have that $ \tilde{a}_{i, \overline{\jmath}(i) } \in \Q \cap U_{i, \overline{\jmath}(i)}$.

\smallskip

SUBCASE 1.2:  $i  \in  X^b$ and  $N_i \neq \emptyset$.

For any $j \in N_i$ choose $\tilde{a}_{i,j} \in U_{i,j} \cap \Q$.

\smallskip
SUBCASE 1.3: $i  \in  X^b$ and  $N_i = \emptyset$.

In this case we have that $a_{i,j} \in \Q$ for any $j \in \{1, \ldots, q-2\}$, so we have to choose no $ \tilde{a}_{i,j}$.

\medskip
CASE 2, i.e.
$\al_{i,q-1}$ is nondegenerate and 
 $\al_{i,q}$ is degenerate.

Analogous to Case 1.

\medskip

CASE 3, i.e.
$\al_{i,q-1}$ and  $\al_{i,q}$ are degenerate.

\smallskip
SUBCASE 3.1: $ i \not \in Y^b \cup Y^c$ and
there does not exist $j
 \in   N_i $ such that $b_j \neq 0 $ and $c_j \neq 0$.

Choose 
$\tilde{a}_{i,j} \in U_{i,j} \cap \Q$ for any $j \in N_i \setminus \{\overline{\jmath}(i), \hat{\jmath}(i)\}$ and define $ \tilde{a}_{i, \overline{\jmath}(i) }$ to be 
\begin{equation} \label{set}
 -\frac{1}{b_{\overline{\jmath}(i)} \wr \tilde{b}_{\overline{\jmath}(i)}} \left[ 
\sum_{j \in  \{1, \ldots, q-2\} \setminus 
\{ \hat{\jmath}(i), \overline{\jmath}(i) \}
}( b_j \wr 
\tilde{b}_j)
\cdot ( a_{i,j} \wr \tilde{a}_{i,j})  -a_{i,q-1} 
\right] .
\end{equation}

 Moreover define 
$ \tilde{a}_{i, \hat{\jmath}(i) }$ to be 
\begin{equation} \label{ot}
 -\frac{1}{c_{\hat{\jmath}(i)} \wr \tilde{c}_{\hat{\jmath}(i)}} \left[ 
\sum_{j \in  \{1, \ldots, q-2\} \setminus 
\{ \hat{\jmath}(i), \overline{\jmath}(i) \}
}( c_j \wr 
\tilde{c}_j)
\cdot ( a_{i,j} \wr \tilde{a}_{i,j})  -a_{i,q} 
\right] .
\end{equation}
By (\ref{ar4}) and (\ref{ar1}), we have that 
$ \tilde{a}_{i, \overline{\jmath}(i) } \in \al_{i, \overline{\jmath}(i)}$ and 
$ \tilde{a}_{i, \hat{\jmath}(i) } \in \al_{i, \hat{\jmath}(i)}$.
\smallskip

SUBCASE 3.2: $ i \not \in Y^b \cup Y^c$ and
there exists $ \overline{\jmath }(i) \in 
  N_i$ such that $b_{\overline{\jmath }(i)} \neq 0$ and 
  $c_{\overline{\jmath }(i)} \neq 0$.

Choose 
$\tilde{a}_{i,j} \in U_{i,j} \cap \Q$ for any $j \in N_i \setminus \{\overline{\jmath}(i)\}$ 
in such way that (\ref{setot}) holds
and define 
$ \tilde{a}_{i, \overline{\jmath}(i) }$ to be one of the members of (\ref{setot}).
Observe that, 
if $ B \cup C
\neq \emptyset$, this is possible by the way we have chosen the $\tilde{b}_j$ and the $\tilde{c}_j$;
if $ B \cup C= \emptyset$, the equation 
(\ref{setot}) has rational coefficients and is solvable,
 because the $ {a}_{i,j}$ for $i \in W$, 
  $j \in N_i \setminus \{\overline{\jmath}(i)\}$
give a solution; so it has a rational solution in 
$ \times_{j \in N_i \setminus \{\overline{\jmath}(i)\}} U_{i,j}$.

\smallskip
SUBCASE 3.3: $ i  \in Y^c \setminus Y^b$
(hence $N_i \neq \emptyset$, $c_j=0 $ $\forall j \in N_i$ and  there exists $\overline{\jmath}(i)
\in N_i$ such that $b_{\overline{\jmath}(i)} \neq 0$).

Choose 
$\tilde{a}_{i,j} \in U_{i,j} \cap \Q$ for any $j \in N_i \setminus \{\overline{\jmath}(i)\}$  
and define $ \tilde{a}_{i, \overline{\jmath}(i) }$ 
as in (\ref{sette}).
\smallskip

SUBCASE 3.4: $ i  \in Y^b \setminus Y^c$
(hence $N_i \neq \emptyset$, $b_j=0 $ $\forall j \in N_i$ and  there exists $\overline{\jmath}(i)
\in N_i$ such that $c_{\overline{\jmath}(i)} \neq 0$).

Analogous to the previous subcase.
\medskip

Finally,
let $H$  be the $ p \times q $ matrix such that, 
 for every $ i=1,\ldots, p$  and $j=1, \ldots, q-2$,
 $$H_{i,j}= \left\{ \begin{array}{ll} 
 \tilde{a}_{i,j} &  \mbox{\rm if } a_{i,j} \in \I  \\ 
 a_{i,j} & \mbox{\rm if } a_{i,j} \in \Q \end{array}
 \right. $$ and such that
$$H^{(q-1)}= 
\sum_{j \in \{1, \ldots, q-2\} |\, b_j \in \Q } b_j H^{(j)}+ \sum_{j  \in \{1, \ldots, q-2\} | \; b_j \in \I}
     \tilde{b}_{j} H^{(j)},$$ 
$$H^{(q)}= 
\sum_{j \in \{1, \ldots, q-2\} |\, c_j \in \Q } c_j H^{(j)}+ \sum_{j  \in \{1, \ldots, q-2\} | \; c_j \in \I}
     \tilde{c}_{j} H^{(j)}.$$ 
     
     By the choice  of $\tilde{b}_j$ for $ j \in B$ and of $\tilde{c}_j$ for $ j \in C$,
     and the choice 
     of $ \tilde{a}_{i,j}$ for $i \in \{k+1, \ldots, p\}
     $, $j \in N_i $, we have 
     that $h_{i,q-1} =a_{i,q-1}$ when $a_{i,q-1} $ is rational, $h_{i,q} =a_{i,q}$ when $a_{i,q} $ is rational, $h_{i,q-1
} \in \al_{i,q-1}$ when $a_{i,q-1} \in \I$ and 
 $h_{i,q
} \in \al_{i,q}$ when $a_{i,q} \in \I$.     
     So the matrix $H$, whose rank is obviously less than or equal to $q-2$, is contained in $ 
     \al \cap M(p \times q, \Q)$.
\end{proof}

\smallskip

From Theorems \ref{Q-R1} and \ref{Q-R111}
 we can easily deduce 
 the following result:

\begin{thm} \label{Q-Rq-1}
Let  $ p \geq q$ and let   
$\al =( [\underline{\alpha}_{i,j}, \overline{\alpha}_{i,j}])_{i,j}$ be a $p \times q$ interval matrix with 
$\underline{\alpha}_{i,j} \leq \overline{\alpha}_{i,j}$ and 
 $\underline{\alpha}_{i,j}, \overline{\alpha}_{i,j} \in \Q$ for any $i,j$. 

(a) 
Suppose there exists $ A \in \al$ with $\rk(A) =q-1$; then 
there exists $ B \in \al \cap M(p \times q, \Q)$ with $\rk(B) =q-1$.

(b) 
Suppose there exists $ A \in \al$ with $\rk(A) =q-2$; then 
there exists $ B \in \al \cap M(p \times q, \Q)$ with $\rk(B) =q-2$.
\end{thm}

\begin{proof}
(a) 
Since there exists $ A \in \al$ with $\rk(A) =q-1$, we can find a $p \times q$ interval matrix $\al'= 
( [\underline{\alpha'}_{i,j}, \overline{\alpha'}_{i,j}])_{i,j}$  with 
$\underline{\alpha'}_{i,j} \leq \overline{\alpha'}_{i,j}$ and 
 $\underline{\alpha'}_{i,j}, \overline{\alpha'}_{i,j} \in \Q$ for any $i,j$ such that $$ A \in \al' \subset \al$$ and 
 \begin{equation} \label{<=q-1}
  \rk(X) \geq q-1 \hspace*{1cm} \forall X \in \al'.
 \end{equation}
 By applying Theorem \ref{Q-R1} to the interval matrix $\al'$, we get that there exists $B \in M(p \times q, \Q) \cap \al'$  with $\rk(B) \leq q-1$; but, by (\ref{<=q-1}), we have that $\rk(B)=q-1$ and we conclude.

(b) Analogous to (a), but we have to use Theorem 
\ref{Q-R111} instead of Theorem \ref{Q-R1}.
\end{proof}

\begin{rem}
We point out that the result above can have some applications: for instance suppose to have a linear subspace $L$ of  dimension $2$ or $q-2$ in  $\R^q$ given as solution set of a linear system  $S$:
$$L=\{x \in \R^q| \;  A x=0\},$$
where $A $ is a $p \times q $ matrix (of rank respectively $q-2$ or $2$ obviously);
  we may want to find a linear subspace $L'$ in
   $\R^q$ with the same dimension as $L$ and given by a linear system $A'x=0$
   such that, for every $i$ and $j$, the entry  $A'_{i,j}$ is equal to $A_{i,j}$ if $A_{i,j}$ is  rational and  the entry $A'_{i,j}$ is rational and in a given interval containing $A_{i,j}$ if $A_{i,j}$ is irrational; the results of this paper allow us to say that this is possible and  the proofs describe also algorithms to find $A'$.
   We observe that the  most ``expensive'' step of the algorithm described in the proof of Theorem \ref{Q-R111} is to calculate the linear system $(X)$ (before Remarks A and B), which requires $O(pq^4) $ elementary operations, and to solve the system $(XN)$,
   which requires $O(p^2q^3) $ elementary operations, while the other steps are less ``expensive'', so we need in all $O(p^2q^3) $ elementary operations.

\end{rem}

\begin{rem} \label{KoBer}
In \cite{Ko},  the authors  exhibited  
   a $12 \times 12 $ sign pattern matrix $Q$ such that 
   there exists a real matrix $B$  
 with $\rk(B)=3$ and sign pattern $Q$   and  
 there does not exist a rational matrix $A$ 
 with $\rk(A) =3$ 
  and  sign pattern $Q$. 
  
  In \cite{Shi} the author
  showed that there exists
   a $p \times q $ sign pattern matrix $Q$ such that there exists 
  a real matrix $B$ 
 with $\rk(B)=q-3$ and sign pattern $Q$    and there does not exist a rational matrix $A$ 
 with $\rk(A) =q-3$ 
  and sign pattern $Q$.

Analogous results are in 
  \cite{Ber}.
  
  The same examples show
  that it is not true for any $r$,
 that, 
 if an interval matrix contains a rank-$r$ real matrix, then it contains   a rank-$r$ rational matrix. In fact
  let $\al$ be an interval matrix containing $B$ of rank $r$ as above and such that, for any $i,j$, we have:

$\al_{i,j}=\{0\}$ if and only if $b_{i,j}=0$,

$\al_{i,j} \subset \R_{>0}$ if and only if $b_{i,j}>0$,

$\al_{i,j} \subset \R_{<0}$ if and only if $b_{i,j}<0$. 

Obviously, since there does not exist a rational matrix with sign pattern $Q$ and rank  $r$,
there does not exist a rational  matrix in $\al $ with rank $r$.
 So Theorem \ref{summarize} is not generalizable to any rank, that is, it is not true for any $r$,
 that, 
 if an interval matrix contains a rank-$r$ real matrix, then it contains   a rank-$r$ rational matrix.

\end{rem}
\bigskip

{\bf Acknowledgments.}
This work was supported by the National Group for Algebraic and Geometric Structures, and their  Applications (GNSAGA-INdAM).


\end{document}